# Entropy of a bit-shift channel

Stan Baggen[1], Vladimir Balakirsky[2], Dee Denteneer[1],
Sebastian Egner[1], Henk Hollmann[1], Ludo Tolhuizen[1] and
Evgeny Verbitskiy[1]

*Philips Research Laboratories Eindhoven and Eindhoven University of Technology*

**Abstract:** We consider a simple transformation (coding) of an iid source called a *bit-shift channel*. This simple transformation occurs naturally in magnetic or optical data storage. The resulting process is not Markov of any order. We discuss methods of computing the entropy of the transformed process, and study some of its properties.

Results presented in this paper originate from the discussions we had at the "Coding Club" – the weekly seminar on coding theory at the Philips Research Laboratories in Eindhoven. Mike Keane, when his active travelling schedule permits, is also attending this seminar. We would like to use this opportunity to thank Mike for his active participation, pleasant and fruitful discussions, his inspiration which we had a pleasure to share.

## 1. Bit-shift channel

In this paper we consider a simplified model for errors occurring in the readout of digital information stored on an optical recording medium like the Compact Disk (CD) or the Digital Versatile Disk (DVD). For more detailed information on optical storage see [9] or [16].

On optical disks the information is stored in a reflectivity pattern. For technical reasons, it is advantageous to use only two states, i.e. "low" and "high" reflectivity. Figure 1 shows the disk surfaces for two types of the DVD's. While the presence of only 2 states greatly simplifies the detection of the state, it reduces the maximum spatial frequency, and hence storage capacity.

In this situation it is better not to encode the information in the reflectivity state itself but rather in the location of the *transitions*: The reflectivity pattern consists of an alternating sequence of "high" and "low" marks of varying length (an integer multiple of some small length unit), while *each* mark exceeds a minimal length, say $d+1$ units. Hence, this "run-length limited" (RLL) encoding makes sure no mark is too short for the disk while the information density is only limited by the accuracy of determining the length of the marks, or equivalently the location of the transitions. For technical reasons (to recover the length unit from the signal itself)

[1]Philips Research Laboratories, Prof. Holstlaan 4, 5656 AA, Eindhoven, The Netherlands, e-mail: `stan.baggen@philips.com` e-mail: `dee.denteneer@philips.com` e-mail: `sebastian.egner@philips.com` e-mail: `henk.d.l.hollmann@philips.com` e-mail: `ludo.tolhuizen@philips.com` e-mail: `evgeny.verbitskiy@philips.com`

[2]TU Eindhoven, 5600 MB, Eindhoven, The Netherlands, e-mail: `V.B.Balakirsky@tue.nl`

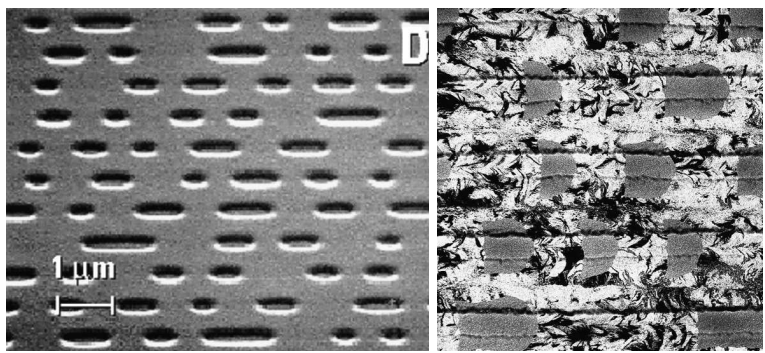

FIG 1. *Images of DVD disks. The left image shows a DVD-ROM. The track is formed by pressing the disk mechanically onto a master disk. On the right is an image of a rewritable disk. The resolution has been increased to demonstrate the irregularities in the track produced by the laser. These irregularities lead to higher probabilities of jitter errors.*

another constraint is imposed: No mark must exceed $k + 1$ units, $k > d$. (For the CD, $(d, k) = (2, 10)$.)

It is customary to describe RLL sequences by their transitions: A $(d, k)$-RLL sequence has at least $d$ and at most $k$ '0's between '1's. So a "high" mark of 4 units, followed by a "low" of 3 units, followed by a "high" of 4 units correspond to the RLL-sequence 100010010001 written to the disk.

At the time the RLL-sequence is read from disk, the transitions (the '1's) might be detected at different positions due to noise, inter-symbol interference, clock jitter, and other distortions. In the simplest version of this "bit-shift channel model" each '1' may be detected one unit early, on time, or one unit late with the probabilities $(\varepsilon, 1 - 2\varepsilon, \varepsilon)$, $0 \leq \varepsilon \leq 1/2$, and the shifts are independent.

More formally, suppose $X$ is the length of a continuous interval of low or high marks on the disk. Then, after reading, the detected length is

$$Y = X + \omega_{left} - \omega_{right}, \tag{1}$$

where $\omega_{left}$, $\omega_{right}$ take values $\{-1, 0, 1\}$. And $\omega = 1, 0, -1$ means that the transition between the "low"-"high" or "high"-"low" runs was detected one time unit too early, correctly, or one unit too late, respectively. Note that for two consecutive intervals $\omega_{right}$ of of the first interval is $\omega_{left}$ of the second. The simplest model for the distribution of time shifts $\omega_{left}$ is to assume that they are independent for different intervals (runs), and

$$\mathbb{P}(\omega_{left} = -1) = \mathbb{P}(\omega_{left} = 1) = \varepsilon, \quad \mathbb{P}(\omega_{left} = 0) = 1 - 2\varepsilon,$$

for some $\varepsilon \in [0, 1/2]$.

An important question then is: Given $(d, k)$, $\varepsilon$, and some distribution for the input sequences (e.g. run-lengths uniformly distributed in $\{d, \ldots, k\}$), what is the mutual information between input and output sequences? In other words, how much can be learned about the input from observing the output, on average. The problem of computing the mutual information is equivalent to computation of the entropy of the output sequence, see [2].

The supremum of this mutual information over all possible measures on space of input sequences is called "channel capacity".

### 1.1. Model

Let us describe the bit-shift channel as a continuous transformation (factor) of a certain *subshift of finite type*. Let $\mathscr{A} = \{d, \ldots, k\}$, where $d, k \in \mathbb{N}$, $d < k$ and $d \geq 2$. The input space then is $\mathscr{A}^{\mathbb{Z}} = \{x = (x_i) : x_i \in \mathscr{A}\}$. Consider also a finite alphabet $\Omega$ with 9 symbols

$$\Omega = \{(-1,-1), (-1,0), (-1,1), (0,-1), (0,0), (0,1), (1,-1), (1,0), (1,1)\}.$$

Finally, consider a subshift of finite type $\Omega_J \subset \Omega^{\mathbb{Z}}$ defined as

$$\Omega_J = \Big\{(\omega_n) \in \Omega^{\mathbb{N}} : \omega_{n,2} = \omega_{n+1,1} \text{ for all } n \in \mathbb{Z}\Big\},$$

where $\omega_n = (\omega_{n,1}, \omega_{n,2})$. The factor map $\phi$ is defined on $\mathscr{A}^{\mathbb{Z}} \times \Omega_J$ as follows: $y = \phi(x, \omega)$ with

$$y_n = x_n + \omega_{n,1} - \omega_{n,2} \quad \text{for all } n. \tag{2}$$

Note that the output space $\mathscr{O} = \phi(\mathscr{A}^{\mathbb{Z}} \times \Omega_J)$ is a subshift of $\mathscr{B}^{\mathbb{Z}}$, with $\mathscr{B} = \{d-2, \ldots, k+2\}$. Clearly, $\mathscr{O} \neq \mathscr{B}^{\mathbb{Z}}$. For example, $(d-2, d-2)$ cannot occur in any output sequence. Indeed, if $y_n = d-2$, then $x_n = d$, $\omega_{n,1} = -1$ and $\omega_{n,2} = 1$. But then, $y_{n+1} = x_{n+1} + \omega_{n+1,1} - \omega_{n+1,2} = x_{n+1} + \omega_{n,2} - \omega_{n+1,2} \geq d + 1 - 1 = d$. With a similar argument, one concludes that for any $L \geq 1$

$$[d-2, \underbrace{d, \ldots, d}_{L \text{ times}}, d-2] \quad \text{or} \quad [d-2, \underbrace{d, \ldots, d}_{L \text{ times}}, d-1]$$

do not occur in any output sequence $y \in \mathscr{O}$. Therefore, there is an infinite number of *minimal* forbidden words, i.e., forbidden words all of whose subwords are allowed. Hence, $\mathscr{O}$ is not a subshift of finite type. However, since $\mathscr{A}^{\mathbb{Z}} \times \Omega_J$ is a subshift of finite type, $\mathscr{O}$ is *sofic* [13].

### 1.2. Capacity of a bit-shift channel

Suppose $\mathbb{J}$ is a measure on $\Omega_J$. For example, in this paper we will be mainly interested in Markov measures $\mathbb{J}_\varepsilon$ on $\Omega_J$, obtained in a natural way from Bernoulli measures on $\{-1, 0, 1\}$ with probabilities $\varepsilon$, $1 - 2\varepsilon$, and $\varepsilon$, respectively. If $\mathbb{P}$ is a translation invariant measure on $\mathscr{A}^{\mathbb{Z}}$, then we obtain a measure $\mathbb{Q}$ on $\mathscr{O}$, which is the push forward of $\mathbb{P} \times \mathbb{J}$. We use a standard notation $\mathbb{Q} = (\mathbb{P} \times \mathbb{J}) \circ \phi^{-1}$.

From the information-theoretical point of view, an important quantity is the *capacity* of the channel. The capacity of a bit-shift channel specified by $\mathbb{J}$ is defined as

$$C_{\text{bitshift}}(\mathbb{J}) = \sup_{\mathbb{P} \in \mathscr{P}(\mathscr{A}^{\mathbb{Z}})} h\Big((\mathbb{P} \times \mathbb{J}) \circ \phi^{-1}\Big) - h(\mathbb{J}), \tag{3}$$

where the supremum is taken over $\mathscr{P}(\mathscr{A}^{\mathbb{Z}})$ – the set of all translation invariant probability measures on $\mathscr{A}^{\mathbb{Z}}$, and $h(\cdot)$ is the entropy.

Even for 'Bernoulli' measures $\mathbb{J}_\varepsilon$, the capacity $C_{\text{bitshift}}(\mathbb{J}_\varepsilon)$ is not known. It is relatively easy to see that the supremum in (3) is achieved. However, the properties of maximizing measures are not known. It is expected that maximizing measures are not Markov of any order. Finally, if one is interested in *topological entropy* of $\mathscr{O}$:

$$h_{top}(\mathscr{O}) = \sup_{\mathbb{Q} \in \mathscr{P}(\mathscr{O})} h(\mathbb{Q}),$$

then $h_{top}(\mathscr{O})$ is easily computable using standard methods [13] or using the efficient numerical approach of [8].

## 2. Entropy of a bit-shift channel

Suppose that $\{X_n\}$ are independent identically distributed random variables taking values in $\mathscr{A} = \{d, \ldots, k\}$, and let $\mathbb{P}$ be the corresponding distribution. What is the entropy of $\mathbb{Q} = (\mathbb{P} \times \mathbb{J}_\varepsilon) \circ \phi^{-1}$?

Note that $(X_n, \omega_n)$ is a Markov chain, and, hence, $Y_n$, given by (2), is a function of a Markov chain.

Let us start by recalling some methods of computing the entropy of processes which are functions of Markov chains. Suppose $\mathscr{X} = \{X_n\}$, $n \in \mathbb{Z}$, is a stationary **ergodic** Markov chain taking values in a finite alphabet $\mathscr{A}$. Let $\phi : \mathscr{A} \to \mathscr{B}$ be some map, and consider a process $\mathscr{Y} = \{Y_n\}$, defined by

$$Y_n = \phi(X_n) \quad \text{for all } n \in \mathbb{Z}.$$

The following result [3, 4], see also [5, Theorem 4.4.1], provides sharp estimates on the entropy of $\mathscr{Y}$.

**Theorem 2.1.** *If $\mathscr{X}$ is a Markov chain and $\mathscr{Y} = \phi(\mathscr{X})$, then for every $n \geq 1$ one has*

$$H(Y_0|Y_1, \ldots, Y_{n-1}, X_n) \leq h(\mathscr{Y}) \leq H(Y_0|Y_1, \ldots, Y_{n-1}, Y_n).$$

*Moreover, as $n \nearrow \infty$*

$$H(Y_0|Y_1, \ldots, Y_{n-1}, X_n) \nearrow h(\mathscr{Y}), \quad H(Y_0|Y_1, \ldots, Y_{n-1}, Y_n) \searrow h(\mathscr{Y}).$$

Birch [3, 4] has shown that under some additional conditions, the convergence is in fact exponential:

$$|h(\mathscr{Y}) - H(Y_0|Y_1, \ldots, Y_{n-1}, X_n)| \leq C\rho^n,$$
$$|h(\mathscr{Y}) - H(Y_0|Y_1, \ldots, Y_{n-1}, Y_n)| \leq C\rho^n,$$

where $\rho \in (0, 1)$ is independent of the factor map $\phi$.

Let us give a proof of Theorem 2.1, since it is very short and provides us with some useful intuition.

*Proof of Theorem 2.1.* An upper estimate of $h(\mathscr{Y})$ in terms of $H(Y_0|Y_1, \ldots, Y_n)$ and the monotonic convergence $H(Y_0|Y_1, \ldots, Y_n)$ to $h(\mathscr{Y})$ are standard facts. For the lower estimate we proceed as follows: for any $m \in \mathbb{N}$ one has

$$\begin{aligned} &H(Y_0|Y_1, \ldots, Y_{n-1}, X_n) \\ &= H(Y_0|Y_1, \ldots, Y_{n-1}, X_n, \ldots, X_{n+m}) &(4)\\ &= H(Y_0|Y_1, \ldots, Y_{n-1}, Y_n, \ldots, Y_{n+m}, X_n, \ldots, X_{n+m}) &(5)\\ &\leq H(Y_0|Y_1, \ldots, Y_{n-1}, Y_n, \ldots, Y_{n+m}), &(6) \end{aligned}$$

where in (4) we used the Markov property of $\mathscr{X}$, and (5), (6) follow from the standard properties of conditional entropies. Since

$$h(\mathscr{Y}) = \lim_{k \to \infty} H(Y_0|Y_1, \ldots, Y_{n-1}, Y_n, \ldots, Y_{n+m}),$$

we obtain the lower estimate of $h(\mathscr{Y})$. Moreover, using standard properties of conditional entropies, we immediately conclude that $H(Y_0|Y_1, \ldots, Y_{n-1}, X_n)$ is monotonically increasing with $n$.

To prove that the lower bound actually converges to $h(\mathscr{Y})$, we proceed as follows. Note that

$$H(X_{n+1}) \geq H(X_{n+1}|Y_n) - H(X_{n+1}|Y_n, \ldots, Y_0)$$
$$= \sum_{i=0}^{n-1} H(X_{n+1}|Y_n, \ldots, Y_{i+1}) - H(X_{n+1}|Y_n, \ldots, Y_i)$$
$$= \sum_{i=0}^{n-1} H(X_{n-i+1}|Y_{n-i}, \ldots, Y_1) - H(X_{n-i+1}|Y_{n-i}, \ldots, Y_0)$$
$$= \sum_{j=1}^{n} H(X_{j+1}|Y_j, \ldots, Y_1) - H(X_{j+1}|Y_j, \ldots, Y_0) = \sum_{j=1}^{n} c_j,$$

where

$$c_j = H(X_{j+1}|Y_j, \ldots, Y_1) - H(X_{j+1}|Y_j, \ldots, Y_0), \quad j = 1, \ldots, n.$$

Since $c_j \geq 0$ and $\sum_{j=1}^{n} c_j < H(X_1) < \infty$ for all $n$, we conclude that $c_n \to 0$ as $n \to \infty$. Moreover,

$$c_n = H(X_{n+1}|Y_n, \ldots, Y_1) - H(X_{n+1}|Y_n, \ldots, Y_0)$$
$$= H(Y_1, \ldots, Y_n, X_{n+1}) - H(Y_1, \ldots, Y_n)$$
$$\quad - H(Y_0, Y_1, \ldots, Y_n, X_{n+1}) + H(Y_0, Y_1, \ldots, Y_n)$$
$$= H(Y_0|Y_1, \ldots, Y_n) - H(Y_0|Y_1, \ldots, Y_n, X_{n+1}).$$

Finally, since $H(Y_0|Y_1, \ldots, Y_n)$ converges to $h(\mathscr{Y})$, so does $H(Y_0|Y_1, \ldots, Y_n, X_{n+1})$. □

Let us conclude this section with one general remark. Suppose $\mathscr{Y}$ is a factor of $\mathscr{X}$, i.e. $\mathscr{Y} = \phi(\mathscr{X})$, where $\mathscr{X}$ is some ergodic process. For $n, m \in \mathbb{N}$, let

$$d_{n,m} = H(Y_0|Y_1, \ldots, Y_n, X_{n+1}, \ldots, X_{n+m}).$$

Note that $d_{n,m} \geq d_{n,m+1} \geq 0$, and hence $\lim_{m \to \infty} d_{n,m} =: D_n$ exists. Since for any $n, m \in \mathbb{N}$

$$d_{n,m} = H(Y_0|Y_1, \ldots, Y_n, X_{n+1}, \ldots, X_{n+m}) \leq H(Y_0|Y_1, \ldots, Y_n, Y_{n+1}, \ldots, Y_{n+m}),$$

we conclude that $D_n \leq h(\mathscr{Y})$. Note also that since $d_{n,m} \leq d_{n+1,m-1}$, one has $D_{n+1} \geq D_n$.

The natural question is under which conditions does $D_n$ converge to $h(\mathscr{Y})$ as $n \to \infty$. For this we need a certain regularity of the conditional probabilities of the $\mathscr{X}$-process. For example, if conditional probabilities are continuous, i.e., if

$$r_n = \sup_{X_0, \ldots, X_n} \sup_{X', X''} \big| \mathbb{P}(X_0|X_1, \ldots, X_n, X'_{n+1}, X'_{n+2}, \ldots)$$
$$- \mathbb{P}(X_0|X_1, \ldots, X_n, X''_{n+1}, X''_{n+2}, \ldots) \big| \to 0, \quad n \to \infty,$$

then $D_n \to h(\mathscr{Y})$. Gibbs measures and $g$-measures (see Section 4) have continuous conditional probabilities.

### 2.1. Entropy via a prefix code

In this section we recall the approach to efficient computation of entropies of factor processes $\mathscr{Y} = \phi(\mathscr{X})$, where $\mathscr{X}$ is Markov, which was originally proposed in [2, 7].

The inequalities of Theorem 2.1 can be rewritten as follows

$$\sum_{y_1^n \in \mathscr{B}^n} \mathbb{P}(y_1^n) H_{\mathbb{P}(\cdot|y_1^n)}(Y_0|X_{n+1}) \leq h(\mathscr{Y}) \leq \sum_{y_1^n \in \mathscr{B}^n} \mathbb{P}(y_1^n) H_{\mathbb{P}(\cdot|y_1^n)}(Y_0), \qquad (7)$$

where we use the following notation

$$y_1^n = (y_1, y_2, \ldots, y_n) \in \mathscr{B}^n,$$
$$\mathbb{P}(y_1^n) = \mathbb{P}(Y_1 = y_1, \ldots, Y_n = y_n),$$
$$\mathbb{P}(\cdot|y_1^n) = \mathbb{P}(\cdot|Y_1 = y_1, \ldots, Y_n = y_n).$$

The subindex $\mathbb{P}(\cdot|y_1^n)$ in (7) stresses that the entropy of $Y_0$ and the conditional entropy of $Y_0$ and $X_{n+1}$ is computed using $\mathbb{P}(\cdot|y_1^n)$.

Note that the sum in (7) is taken over elements of a partition of $\mathscr{B}^{\mathbb{N}}$ into cylinders of length $n$:

$$\mathscr{U}_n = \left\{[y_1^n] \,\middle|\, y_1^n \in \mathscr{B}^n\right\}, \quad [y_1^n] = \{\tilde{y} \in \mathscr{B}^{\mathbb{Z}} : \tilde{y}_1 = y_1, \ldots, \tilde{y}_n = y_n\}.$$

In fact, an estimate similar to (7) holds for **any** partition of $\mathscr{B}^{\mathbb{Z}}$ into cylindric sets, see [7, Theorem 1].

**Theorem 2.2.** *Let $\mathscr{W}$ be a finite partition of $\mathscr{B}^{\mathbb{Z}}$ into cylindric sets:*

$$\mathscr{W} = \Big\{[\mathbf{w}_i], \ \mathbf{w}_i = (w_{i,1}, \ldots, w_{i,l_i})\Big\}_{i=1}^M.$$

*Then*

$$\sum_{\mathbf{w} \in \mathscr{W}} h_1(\mathbf{w}) \leq h(\mathscr{Y}) \leq \sum_{\mathbf{w} \in \mathscr{W}} h(\mathbf{w}), \qquad (8)$$

*where*

$$h_1(\mathbf{w}) = \mathbb{P}(Y_1 \ldots Y_{|\mathbf{w}|} = \mathbf{w}) H(Y_0 | Y_1 \ldots Y_{|\mathbf{w}|} = \mathbf{w}, X_{|\mathbf{w}|+1}),$$
$$h(\mathbf{w}) = \mathbb{P}(Y_1 \ldots Y_{|\mathbf{w}|} = \mathbf{w}) H(Y_0 | Y_1 \ldots Y_{|\mathbf{w}|} = \mathbf{w}).$$

Theorem (2.2) leads to the following algorithm. Suppose $\mathscr{W}$ is some partition into cylinders. We can refine the partition $\mathscr{W}$ by removing a certain word $\mathbf{w}$ from $\mathscr{W}$ and adding all words of the form $\mathbf{w}b$, where $b \in \mathscr{B}$, i.e.,

$$\mathscr{W}' = \mathscr{W} \setminus \{\mathbf{w}\} \cup \{\mathbf{w}b | \, b \in \mathscr{B}\}. \qquad (9)$$

Suppose $\{\mathscr{W}_k\}_{k \geq 1}$ is a sequence of partitions such that for each $k$, $\mathscr{W}_{k+1}$ is a refinement of $\mathscr{W}_k$ as in (9), and at each step a word $\mathbf{w} \in \mathscr{W}_k$ is selected such that

$$h(\mathbf{w}) - h_1(\mathbf{w}) = \max_{u \in \mathscr{W}_k} \Big(h(\mathbf{u}) - h_1(\mathbf{u})\Big). \qquad (10)$$

The greedy strategy (10), as well as some other strategies (e.g, *uniform*, $|\mathbf{w}| = \min_{\mathbf{u} \in \mathscr{W}} |\mathbf{u}|$), guarantees the convergence of the upper and lower estimates in (8), i.e.,

$$\lim_{k \to \infty} \sum_{\mathbf{w} \in \mathscr{W}_k} \Big(h(\mathbf{w}) - h_1(\mathbf{w})\Big) = 0.$$

## 2.2. Entropy via renewal times.

As before, suppose that $\{X_n\}$ are independent and identically distributed in $\{d, \ldots, k\}$ with $\mathbb{P}(X_i = \ell) = p_\ell$, $\ell = d, \ldots, k$. Assume also that $p_d > 0$.

Another method for estimating the entropy is based on the following observation. Suppose $Y_n = d - 2$ for some $n$. This implies that $X_n = d$, $\omega_n = (-1, 1)$. Since the sequence $\{\omega_k\}$ forms a Markov chain, $(\ldots, \omega_1, \ldots, \omega_{n-1})$ and $(\omega_{n+1}, \omega_{n+2}, \ldots)$, are independent given $\omega_n$. Therefore, since $Y_n = d-2$ implies $\omega_n = (-1, 1)$, we conclude that $(\ldots, \omega_1, \ldots, \omega_{n-1})$ and $(\omega_{n+1}, \omega_{n+2}, \ldots)$ are independent given $Y_n = d - 2$. Moreover, since $X_n$ form an iid sequence, $(\ldots, Y_{n-2}, Y_{n-1})$ and $(Y_{n+1}, Y_{n+2}, \ldots)$ are also independent given $Y_n = d - 2$.

Consider our subshift $\mathscr{O}$, and a set $C = [d-2] = \{y \in \mathscr{O} : y_0 = d - 2\}$. Let $S : \mathscr{O} \to \mathscr{O}$ be a left shift, and consider an induced map $S_C$ on $C$:

$$S_C(y) = S^{R_C(y)}(y),$$

where $R_C(y) = \min\{k \geq 1 : y_k = d - 2\}$. On $C$, the induced map $S_C$ has a natural Bernoulli partition

$$\Big\{[d-2, y_1, \ldots, y_r, d-2] : \ y_j \in \mathscr{B}, \ y_j \neq d-2, \ j = 1, \ldots, r, \ r \in \mathbb{N}\Big\}.$$

Finally, by the Abramov formula [1]

$$h(\mathbb{Q}) = -\sum_{r=1}^{\infty} \sum_{y_1, \ldots, y_r \neq d-2} \mathbb{Q}([d-2, y_1, \ldots, y_r, d-2]) \log \mathbb{Q}\,[d-2, y_1, \ldots, y_r, d-2])$$
$$+ \mathbb{Q}([d-2]) \log \mathbb{Q}([d-2]). \tag{11}$$

Computation of entropy of images of Markov measures using the renewal times and induced map was used in the past, see e.g. [15]. However, in the case of bit-shift channel, the method based on (11) is extremely inefficient.

## 3. Numerics

For illustration we present a numerical computation of the entropy using the prefix code method described in Section 2.1.

The algorithm constructs a sequence of refined partitions $\mathscr{W}_k$ as described above. A particularly useful strategy is given by (10). This "greedy" heuristics selects the cylinder most responsible for the difference in upper and lower bound, in the hope that refining this cylinder will tighten the bounds quickly. This strategy is not optimal (as can be shown by example) but it has three advantages. Firstly, the bounds converge (eventually). Secondly, if in a particular word $\mathbf{w} \in \mathscr{W}$, the last symbol is the "renewal" symbol $d-2$ (similarly $k+2$), this word will never be refined again. Thirdly, the next cylinder to expand can be found quickly by representing $\mathscr{W}$ as a "priority queue" data structure.

For illustration, we run the algorithm for the model of the jitter channel described in Section 1.1. The parameters are inspired by the Compact Disc: The error-correction and modulation system of the CD essentially produces an RLL-sequence with parameters $(d, k) = (2, 10)$. We model the run-lengths as independent identically distributed random variables with probabilities $p_\ell = p_2 \gamma^{\ell-2}$, $\ell \in \{2, \ldots, 10\}$, where $\gamma = 0.658$ and $p_2$ is chosen such that $\sum p_\ell = 1$. This truncated geometric model with $\gamma = 0.658$ is a very good approximation of the (marginal) run-length distribution observed on the CD.

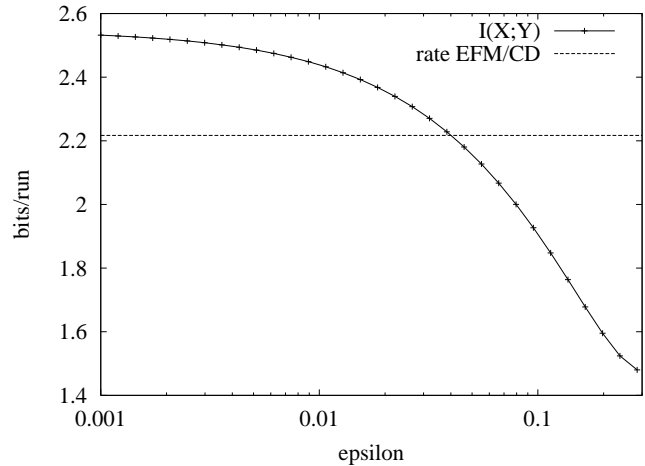

FIG 2. *Mutual information $I(\mathcal{Y}; \mathcal{X}) = h(\mathcal{Y}) - h(\mathbb{J}_\varepsilon)$ as a function of $\varepsilon$, for $(d, k) = (2, 10)$ and the truncated geometric distribution for $\mathcal{X}$.*

Figure 2 shows the mutual information

$$I(\mathcal{Y}; \mathcal{X}) = h(\mathcal{Y}) - h(\mathbb{J}_\varepsilon) = h(\mathcal{Y}) + 2\varepsilon \log \varepsilon + (1 - 2\varepsilon) \log(1 - 2\varepsilon)$$

as a function of $\varepsilon$. The horizontal line represents the rate designed for the last stage of the encoding used in the CD (the so-called EFM code). If the jitter is so strong that the mutual information drops below this rate, reliable decoding is impossible. In practice, similar plots are used to evaluate the performance of particular encoding schemes with respect to various distortions introduced by the physical channel.

Figure 3 compares the greedy and uniform heuristics. The standard estimate $H(Y_0|Y_1, \ldots, Y_n)$ in fact corresponds to the uniform refinement. Observe a superior rate of convergence for the greedy refinement strategy.

## 4. Thermodynamics of jittered measures

Bernoulli and Markov measures belong to a wider class of the so-called Gibbs measures. Bernoulli and Markov measures are also examples of $g$-measures.

In the seminal paper [10] M. Keane introduced a class of *g-measures*. These are the measures whose conditional probabilities are given by a continuous and strictly positive function $g$. For subshifts of finite type, the theory of $g$-measures is extensive. For sofic subshifts, the problem of defining $g$-measures is much more complicated. For the first results see the paper by W. Krieger [11] in this volume.

The thermodynamic formalism allows to look at Gibbs measures from two different sides. First of all, *locally*, through the conditional probabilities; and secondly, *globally*, through the variational principles.

Contrary to the class of $g$-measures, the class of Gibbs measures for a sofic subshift is well defined. The natural question is whether a "jittered" measure $\mathbb{Q} = (\mathbb{P} \times \mathbb{J}_\varepsilon) \circ \phi^{-1}$ is Gibbs. If the measure is Gibbs and the potential is identified, then,

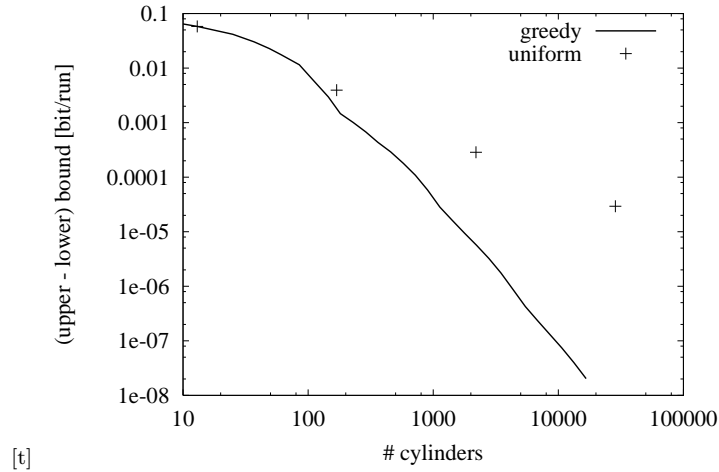

FIG 3. *Difference between upper and lower bounds on entropy as a function of $|\mathcal{W}|$, the number of cylinders, in partitions built by the greedy and uniform refinement strategies.*

using the variational principle, we obtain another method of computing the entropy of $\mathbb{Q}$.

The subshift $\mathscr{O} \subset \mathscr{B}^{\mathbb{Z}}$ satisfies a *specification* property (as a factor of a subshift of finite type $\mathscr{A}^{\mathbb{Z}} \times \Omega_J$ which has a specification property [6]). Hence the results of [17] on existence of Gibbs measures for expansive dynamical systems with the specification property are applicable. If $\mathbb{Q}$ would be a Gibbs measure for potential $f$ from the Bowen class $\mathcal{V}(\mathscr{O})$, then there would exist positive constants $c$, $C$ such that for any $n > 0$ and every $y \in \mathscr{O}$

$$c \leq \frac{\mathbb{Q}([y_0, \ldots, y_n])}{\exp(\sum_{k=0}^{n} f(S^i(y)) - (n+1)P(f))} \leq C, \qquad (12)$$

where $S : \mathscr{O} \to \mathscr{O}$ is a left shift and $P(f)$ is the topological pressure of $f$. As a corollary of (12) one easily concludes that

$$\log \mathbb{Q}(y_0 | y_1, \ldots, y_n) = \log \frac{\mathbb{Q}([y_0, \ldots, y_n])}{\mathbb{Q}([y_1, \ldots, y_n])}$$

should be bounded, which is not the case, see [18]. Hence, $\mathbb{Q}$ is not Gibbs for any potential from the large class of potentials $\mathcal{V}(\mathscr{O})$. Examples of measures $\mathbb{Q}$ such that estimates similar to (12) hold for some continuous $f$ and subexponential bounds $c_n$ and $C_n$ ($\lim_n n^{-1} \log c_n = \lim_n n^{-1} \log C_n = 0$) have been considered [14, 20], and were shown to be *weakly Gibbs*. It is not known whether $\mathbb{Q}$ for the bit-shift channel is weakly Gibbs for some continuous potential $f$.

Nevertheless, the thermodynamic formalism could be useful in estimating the capacity of the bit-shift channel. We recall the notion of a compensation function and some results summarized in [19]. First of all, we define the topological pressure of real-valued continuous functions defined on $\mathscr{A}^{\mathbb{Z}} \times \Omega_J$ and $\mathscr{O}$. If $f \in C(\mathscr{A}^{\mathbb{Z}} \times \Omega_J)$, $g \in C(\mathscr{O})$, the topological pressures of $f$ and $g$ are defined as

$$P(f | \mathscr{A}^{\mathbb{Z}} \times \Omega_J) = \sup_{\mathbb{S}} \left( h(\mathbb{S}) + \int_{\mathscr{A}^{\mathbb{Z}} \times \Omega_J} f \, d\mathbb{S} \right), \quad P(g | \mathscr{O}) = \sup_{\mathbb{Q}} \left( h(\mathbb{Q}) + \int_{\mathscr{O}} g \, d\mathbb{Q} \right),$$

where the suprema are taken over all translation invariant measures on $\mathscr{A}^{\mathbb{Z}} \times \Omega_J$ and $\mathscr{O}$, respectively. A measure $\mathbb{S}$ on $\mathscr{A}^{\mathbb{Z}} \times \Omega_J$ is called an *equilibrium state* for $f \in C(\mathscr{A}^{\mathbb{Z}} \times \Omega_J)$ if

$$P(f|\mathscr{A}^{\mathbb{Z}} \times \Omega_J) = h(\mathbb{S}) + \int f d\mathbb{S}. \tag{13}$$

We define equilibrium states on $\mathscr{O}$ in a similar way. It is well known that every measure is an equilibrium state: for every translation invariant measure $\mathbb{S}$ on $\mathscr{A}^{\mathbb{Z}} \times \Omega_J$ one can find a continuous function $f : \mathscr{A}^{\mathbb{Z}} \times \Omega_J \to \mathbb{R}$ such that (13) holds. Moreover, for any $\mathbb{S} = \mathbb{P} \times \mathbb{J}_\varepsilon$, such an $f$ is of a special form

$$f(x, \omega) = \tilde{f}(x) + j_\varepsilon(\omega),$$

where $\tilde{f} : \mathscr{A}^{\mathbb{Z}} \to \mathbb{R}$ and $j_\varepsilon : \Omega_J \to \mathbb{R}$ are continuous functions. (In fact, $j_\varepsilon$ can be found explicitly.)

A continuous function $F : \mathscr{A}^{\mathbb{Z}} \times \Omega_J \to \mathbb{R}$ is a *compensation* function if

$$P(F + g \circ \phi|\mathscr{A}^{\mathbb{Z}} \times \Omega_J) = P(g|\mathscr{O})$$

for all $g \in C(\mathscr{O})$. Compensation functions exist for factor maps defined on shifts of finite type [19].

An important result is the so-called *relative variational principle* [12, 19], which in our notation states that $F$ is a compensation function if and only if for any invariant measure $\mathbb{Q}$ on $\mathscr{O}$ one has

$$h(\mathbb{Q}) = \sup_{\mathbb{S}} \Big( h(\mathbb{S}) + \int F \, d\mathbb{S} \,\Big|\, \mathbb{S} \circ \phi^{-1} = \mathbb{Q} \Big).$$

Suppose $F$ is a compensation function, then for $\mathbb{Q} = (\mathbb{P} \times \mathbb{J}_\varepsilon) \circ \phi^{-1}$ we obtain

$$h(\mathbb{Q}) \geq h(\mathbb{P} \times \mathbb{J}_\varepsilon) + \int F \, d(\mathbb{P} \times \mathbb{J}_\varepsilon) = h(\mathbb{P}) + h(\mathbb{J}_\varepsilon) + \int F_\varepsilon d\mathbb{P}, \tag{14}$$

where $F_\varepsilon(x) = \int_{\Omega_J} F(x, \omega) \mathbb{J}_\varepsilon(d\omega)$. For the capacity of the bit-shift channel we obtain the following lower estimate

$$C_{\text{bitshift}}(\mathbb{J}_\varepsilon) = \sup_{\mathbb{Q} = (\mathbb{P} \times \mathbb{J}_\varepsilon) \circ \phi^{-1}} h(\mathbb{Q}) - h(\mathbb{J}_\varepsilon)$$

$$\geq \sup_{\mathbb{P}} \Big( h(\mathbb{P}) + \int F_\varepsilon d\mathbb{P} \Big) = P(F_\varepsilon|\mathscr{A}^{\mathbb{Z}}). \tag{15}$$

An interesting question is whether the inequalities in (14) and (15) are strict. The inequality (14) is most probably strict in the generic situation. Indeed, by Corollary 3.4 [19], if $\mathbb{Q}$ is an equilibrium state for $g$ on $\mathscr{O}$, and $\mathbb{S}$ is such that $\mathbb{S} \circ \phi^{-1} = \mathbb{Q}$ and

$$h(\mathbb{Q}) = h(\mathbb{S}) + \int F \, d\mathbb{S},$$

then $\mathbb{S}$ is an equilibrium state for $F + g \circ \phi$, and conversely. On the other hand if $\mathbb{S} = \mathbb{P} \times \mathbb{J}_\varepsilon$, then $\mathbb{S}$ is an equilibrium state for $f(x, \omega) = \tilde{f}(x) + j_\varepsilon(\omega)$. Therefore, for the equality in (14), it is necessary that $F(x, \omega) + (g \circ \phi)(x, \omega)$ and $\tilde{f}(x) + j_\varepsilon(\omega)$ are *physically equivalent*, i.e., have the same set of equilibrium states. In fact, it is quite difficult to imagine how for a given compensation function $F$ of the bit-shift channel and a generic $g$ one could find $\tilde{f}$ to satisfy the requirement of physical equivalence. On the other hand, it is not very difficult to see that in fact

$$C_{\text{bitshift}}(\mathbb{J}_\varepsilon) = \sup_F P(F_\varepsilon|\mathscr{A}^{\mathbb{Z}}), \tag{16}$$

where the supremum is taken over all compensation functions $F$. Indeed, suppose $\mathbb{Q}^* = (\mathbb{P}^* \times \mathbb{J}_\varepsilon) \circ \phi^{-1}$ is a 'maximal' ergodic measure, i.e.,

$$C_{\text{bitshift}}(\mathbb{J}_\varepsilon) = h(\mathbb{Q}^*) - h(\mathbb{J}_\varepsilon).$$

Then there exist continuous functions $g^* \in C(\mathscr{O})$ and $\tilde{f}^* \in C(\mathscr{A}^{\mathbb{Z}})$ such that $\mathbb{Q}^*$ and $\mathbb{P}^*$ are equilibrium states for $g^*$ and $\tilde{f}^*$, respectively. But then

$$F(x,\omega) = \tilde{f}^*(x) + j_\varepsilon(\omega) - (g^* \circ \phi)(x,\omega)$$

is the compensation function for which the maximum in (16) is attained. Thus methods for dealing with factor systems developed in dynamical systems, could be applied to estimate channel capacities. The practicality of such estimates depends strongly on whether one is able to understand the structure of a class of compensation function for a given channel. Probably, in many concrete cases, a relatively large family of compensation functions will suffice as well.

## References

[1] ABRAMOV, L. M. (1959). The entropy of a derived automorphism. *Dokl. Akad. Nauk SSSR* **128**, 647–650. MR0113984

[2] BAGGEN, S., AND BALAKIRSKY, V. (2003). An efficient algorithm for computing the entropy of output sequences for bitshift channels. *Proc. 24$^{th}$ Int. Symposium on Information Theory in Benelux*, 157–164.

[3] BIRCH, J. J. (1962) Approximation for the entropy for functions of Markov chains. *Ann. Math. Statist.* **33**, 930–938. MR0141162

[4] BIRCH, J. J. (1963) On information rates for finite-state channels. *Information and Control* **6**, 372–380. MR0162651

[5] COVER, T. M., AND THOMAS, J. A. (1991). *Elements of Information Theory*. Wiley Series in Telecommunications. John Wiley & Sons Inc., New York. MR1122806

[6] DENKER, M., GRILLENBERGER, C., AND SIGMUND, K. (1976). *Ergodic Theory on Compact Spaces*. Springer-Verlag, Berlin. Lecture Notes in Mathematics **527**. MR0457675

[7] EGNER, S., BALAKIRSKY, V., TOLHUIZEN, L., BAGGEN, S., AND HOLLMANN, H. (2004) On the entropy rate of a hidden Markov model. *Proceedings International Symposium Information Theory, ISIT 2004*.

[8] FROYLAND, G., JUNGE, O., AND OCHS, G. (2001). Rigorous computation of topological entropy with respect to a finite partition. *Phys. D* **154**, 1–2, 68–84. MR1840806

[9] IMMINK, K. (1999). *Codes for Mass Data Storage Systems*. Shannon Foundation EEE Publishers, The Netherlands.

[10] KEANE, M. (1972). Strongly mixing $g$-measures. *Invent. Math.* **16**, 309–324. MR0310193

[11] KRIEGER, W. (2006) On $g$-functions for subshifts. In *Dynamics and Stochastics*, IMS Lecture Notes-Monograph Series, Vol. 48, 306–316.

*Entropy of a bit-shift channel* 285[12] LEDRAPPIER, F., AND WALTERS, P. (1977). A relativised variational principle for continuous transformations. *J. London Math. Soc. (2)* **16**, 3, 568–576. MR0476995

[13] LIND, D., AND MARCUS, B. (1995). *An Introduction to Symbolic Dynamics and Coding.* Cambridge University Press. MR1369092

[14] MAES, C., REDIG, F., TAKENS, F., VAN MOFFAERT, A., AND VERBITSKI, E. (2000). Intermittency and weak Gibbs states. *Nonlinearity* **13**, 5, 1681–1698. MR1781814

[15] MARCUS, B., PETERSEN, K., AND WILLIAMS, S. ( 1984). Transmission rates and factors of Markov chains. In *Conference in Modern Analysis and Probability (New Haven, Conn., 1982)*, vol. 26 of *Contemp. Math.* Amer. Math. Soc., Providence, RI, pp. 279–293. MR1369092

[16] MARCUS, B. H., ROTH, R. M., AND SIEGEL, P. H. (2001). *An Introduction to Coding of Constrained Systems.* Lecture Notes, fifth edition.

[17] RUELLE, D. (1992). Thermodynamic formalism for maps satisfying positive expansiveness and specification. *Nonlinearity* **5**, 6, 1223–1236. MR1192516

[18] VAN ENTER, A. C. D., AND VERBITSKIY, E. A. (2004). On the variational principle for generalized Gibbs measures. *Markov Process. Related Fields* **10**, 3, 411–434. MR2097865

[19] WALTERS, P.(1986). Relative pressure, relative equilibrium states, compensation functions and many-to-one codes between subshifts. *Trans. Amer. Math. Soc.* **296**, 1, 1–31. MR837796

[20] YURI, M. (1999). Thermodynamic formalism for certain nonhyperbolic maps. *Ergodic Theory Dynam. Systems* **19**, 5, 1365–1378. MR1721626